\documentclass[12pt]{amsart}
\usepackage{amsmath,amsthm,amssymb}
\usepackage{fullpage}
\usepackage{amssymb}
\usepackage{pstricks}
\usepackage{graphics}
\usepackage{setspace}
\usepackage{multicol}
\usepackage[centertags]{amsmath}
\usepackage{amsfonts,stmaryrd,yfonts}

\newtheorem{teo}{Theorem}[section]

\newtheorem{coro}[teo]{Corollary}
\newtheorem{lema}[teo]{Lemma}
\newtheorem{ejem}[teo]{Example}
\newtheorem{ejems}[teo]{Examples}

\newtheorem{rem}[teo]{Remark}
\newtheorem{prop}[teo]{Proposition}
\newtheorem{defi}[teo]{Definition}

\DeclareMathOperator{\lex}{lex}
\DeclareMathOperator{\uf}{\mathcal F}
\DeclareMathOperator{\finite}{\mathbb{N}^{^{[<\infty]}}}

\begin{document}
\hyphenation{hy-po-the-sis exist}

\title{On topological properties  of families of finite sets}
\author{Claribet Pi\~{n}a}
\address{Departamento de Matem\'aticas, Facultad de
Ciencias,Universidad de Los Andes, M\'erida, 5101, Venezuela}
\email{claribet@ula.ve}

\author{Carlos Uzc\'ategui}\thanks{We thank the partial support provided by the  University of Los Andes CDCHTA grant ???????}
\address{Departamento de Matem\'aticas, Facultad de
Ciencias, Universidad de Los Andes, M\'erida, 5101, Venezuela}
\email{uzca@ula.ve}
\date{}



\begin{abstract}
We present results about the Cantor-Bendixson index of some
subspaces of a uniform family $\mathcal F$ of finite subsets of
natural numbers with respect to the lexicographic order topology. As
a corollary of our results we get that for any $\omega$-uniform
family $\mathcal F$ the restriction ${\mathcal F}\upharpoonright M$
is homeomorphic to $\uf$ iff $M$ contains intervals of arbitrary
length of consecutive integers. We show the connection of these
results with a topological partition problem of uniform families.
\end{abstract}

\maketitle

\section{Introduction}

A partition problem for topological spaces is as follows: Given
spaces $X$ and $Y$ and a partition of $X$ into two pieces, is
there a topological copy of $Y$ inside one of the pieces? When the
answer is positive, it is denoted by $X\rightarrow(Y)^1_2$ (see
\cite{Komjath_weiss1987,Weiss90} for more information about this
type of problems). We will be interested in the case $X=Y$. A
result of Baumgartner \cite{Baumgart86} solves this partition
problem  when $X$ is a countable ordinal space $\alpha$. Namely,
he showed that for a countable ordinal $\alpha$,
$\alpha\rightarrow (\alpha)^1_2$ iff $\alpha$ is of the form
$\omega^{\omega^\beta}$.

Any countable ordinal is the order type of a uniform family
$\mathcal F$ of finite subsets of natural numbers
lexicographically ordered. A typical uniform family of order type
$\omega^k$ is the collection of $k$-elements subsets of $\mathbb
N$. Thus a partition of a countable ordinal space can be regarded
as a partition of a uniform family endowed with the lexicographic
order topology (the relevant definitions  are given on section
\ref{terminology}).

Families of finite sets has been the focus of Ramsey theory for a
long time \cite{Todor2010}. A well known result of Nash-Williams
says that for any uniform family $\uf$ on $\mathbb N$ and any subset
$\mathcal B$ of $\mathcal F$ there is an infinite set $A\subseteq
\mathbb N$ such that either ${\mathcal F}\upharpoonright A\subseteq
{\mathcal B}$ or ${\mathcal F}\upharpoonright A\cap {\mathcal
B}=\emptyset$ (see \cite{Todor2010}) where ${\mathcal
F}\upharpoonright A$ is the collection of elements of $\mathcal F$
that are subsets of $A$. This theorem solves the topological
partition problem for $\mathcal F$, if the topological type of
${\mathcal F}\upharpoonright A$ and $\mathcal F$ are the same. This
was the starting point for this research. We soon realized that
${\mathcal F}\upharpoonright A$ could be a discrete subspace of
$\mathcal F$ and hence Baumgartner's theorem is not a corollary of
the Nash-Williams's theorem. In fact, given a uniform family $\uf$,
there is ${\mathcal B}\subset \uf$ such that ${\mathcal
F}\upharpoonright A$ is a discrete subset of $\uf$ for every set $A$
homogeneous for the partition given by $\mathcal B$ (i.e. for any
$A$ satisfying the conclusion of Nash-Williams's theorem applied to
$\uf$ and $\mathcal B$) (see Example \ref{nash-williams}).
Nevertheless, it is natural to wonder about the topological type of
${\mathcal F}\upharpoonright A$. The objective of this paper is to
present an analysis of the Cantor-Bendixson index of ${\mathcal
F}\upharpoonright A$ as a subspace of a uniform family $\mathcal F$.
Notice that ${\mathcal F}\upharpoonright A$ has the same order type
of $\mathcal F$, but the topological type varies considerably
depending on the set $A$. Hence the difficulty lies on the fact that
we are using on ${\mathcal F}\upharpoonright A$ the subspace
topology.

To give an example of the results presented in this paper, we
recall a typical $\omega$-uniform family, the so called {\em
Schreier barrier}:
\[
\mathcal{S}=\{t\in \finite:\; |t|=\min(t)+1\}.
\]
It is known that $\mathcal S$ is homeomorphic to $\omega^\omega$.
We will show that ${\mathcal S}\upharpoonright M$ contains a
topological copy of  $\mathcal S$ iff $M$ contains arbitrarily
long intervals of consecutive natural number. In fact, this result holds for any $\omega$-uniform family. 

We show a partial generalization of the previous result for $\alpha$-uniform families. However, it is still open the general problem of characterizing (by combinatorial means) the class of infinite subsets $M$ of $\mathbb N$ for which ${\mathcal F}\upharpoonright M$ contains a topological copy of $\uf$ when $\uf$ is $\alpha$-uniform on $\mathbb N$ for some $\alpha < \omega_1$.

Finally, we mention that uniform families have been extensively used in the theory of Banach spaces (see, e.g., \cite{argyros2005ramsey, LopezManoussakis2008}).

The paper is organized as follows. In section 2 we introduce the
terminology and some preliminary facts. In section 3 we study the
Cantor-Bendixson derivatives of uniform families. In section 4 we
introduce the type of sets $M$ such that the restriction
$\uf\upharpoonright M$ has the same Cantor-Bendixson index as
$\uf$.  Finally, in section 5 we present the main results about
when $\uf\upharpoonright M$ contains a topological copy of $\uf$.

\section{Preliminaries}
\label{terminology}
We denote by $\finite$ the collection of all
finite subsets of $\mathbb N$. If $M$ is a set, $M^{[k]}$ denotes
the collection of all $k$-elements subsets of $M$. By
$M^{[\infty]}$ we denote the collection of all infinite subsets of
$M$.

The lexicographic order $<_{\lex}$ over $\finite$ is defined as
follows: Given $s,t\in \mathbb{N}^{[<\infty]}$ we put $s<_{\lex}t$
iff $\min(s{\scriptstyle \triangle} t) \in s$.

We write $s\sqsubseteq t$ when there is $n\in \mathbb N$ such that
$s=t\cap \{0,1,\cdots, n\}$ and we say that $s$ is an initial
segment of $t$.  A collection $\mathcal F$ of finite subsets of
$\mathbb N$ is a {\em front} on $M$ if satisfies the following
conditions: (i) Every two elements of $\mathcal F$ are
$\sqsubseteq$-incomparable. (ii) Every infinite subset $N$ of $M$
has an initial segment in $\mathcal F$.

Given $\mathcal{F}\subseteq \finite$ and  $u\in\finite$, let
\[
\mathcal{F}_{u}=\{s\in\finite:\; u\cup s\in  \mathcal{F}, \;
\max(u)<\min(s)\}.
\]
For convenience, we set $\max(\emptyset)=-1$; in particular,
$\uf_\emptyset=\uf$.

For $M$  an infinite subset of $\mathbb N$, let
\[
\uf\upharpoonright M=\{s\in \uf:\;s\subset M\}.
\]
We put $M/k=\{n\in M:\; k<n\}$. If $u$ is a finite set and
$n=\max(u)$, we put $M/u=M/n$. 

The notion of an $\alpha$-uniform
family on an infinite set $M$  was introduced by P. Pudl\'{a}k and V. R\"{o}dl \cite{PudlakRodl1982}. It is is defined by recursion.

\begin{itemize}
\item[(i)] $\{\emptyset\}$ is the unique $0$-uniform family on
$M$.

\item[(ii)] ${\mathcal F}\subseteq \finite$ is said to be
$(\alpha+1)$-uniform on $M$, if $\mathcal{F}_{\{n\}}$ is
$\alpha$-uniform on $M/n$ for all $n\in M$.

\item[(iii)] If $\alpha$ is a limit ordinal, we say that $\mathcal
F$ is $\alpha$-uniform on $M$, if there is an increasing sequence
$(\alpha_k)_{k\in M}$ converging to $\alpha$ such
that $\mathcal{F}_{\{k\}}$ is $\alpha_k$-uniform on $M/k$ for all
$k\in M$.
\end{itemize}

For $k\in\mathbb N$, $M^{[k]}$ is the unique $k$-uniform family on
$M$. The following collection is an $\omega$-uniform family on
$\mathbb N$, called {\em Schreier barrier}:
\[
\mathcal{S}=\{t\in \finite:\; |t|=\min(t)+1\}.
\]
Collection of finite sets similar to the Schreier barrier were studied in  \cite{Frmaki2002}.

We say that $\uf$ is uniform on $M$ when it is $\alpha$-uniform on
$M$ for some $\alpha$. Notice that if $\uf$ is uniform on $M$,
then $\uf_u$ es uniform on $M/u$. 

The following result is well known \cite{argyros2005ramsey}.

\begin{teo}
Let $\uf$ be an $\alpha$-uniform family over $M$. Then $\uf$ is a
front over $M$ and $\uf\upharpoonright N$ is $\alpha$-uniform
 over $N$ for all infinite $N\subseteq M$.
\end{teo}

Given a front  $\mathcal{F}$ on a final segment $S$ of $\mathbb
N$. For $n\in S$,  we denote by $t_n^{\mathcal F}$ the unique
element of $\mathcal{F}$ verifying
\[
t_n^{\mathcal F}\sqsubseteq \{n,n+1,n+2,\dots\}.
\]
In the sequel, the sets $t_n^{{\mathcal F}_u}$ will be very
useful. In particular, we remark that given a finite set $u\subset
S$ and $n\in  S/u$, there is a unique $m$ such that
\[
u\cup t_n^{{\mathcal F}_u}=u\cup\{n,n+1,\dots,n+m\}\in {\mathcal F}.
\]

Notice that if $s\in \mathcal F$ and $n=\min(s)$, then
\[
t_n^{\mathcal F}\leq_{\lex}s <_{\lex} t_{n+1}^{\mathcal F}.
\]

Given two families $\mathcal{F}$ and $\mathcal{G}$ of finite sets,
define $\mathcal{F} \varoplus \mathcal{G}$ as follows:
\[
\mathcal{F} \varoplus \mathcal{G}=\{s\cup t:\; s\in \mathcal{G},\;
t\in \mathcal{F}\;\mbox{and $\max(s)<\min(t)$}\}.
\]
If $\uf$ is $\alpha$-uniform and $\mathcal G$ is $\beta$-uniform,
then $\mathcal{F} \varoplus \mathcal{G}$ is
$(\alpha+\beta)$-uniform. Notice that if $\uf$ is a front over a
final segment $S$ of $\mathbb N$, then $t_n^{\mathcal
F}=\min(\mathcal{F}_{\{n\}}\varoplus \{\{n\}\},<_{\lex})$ for all
$n\in S$.

The following result is well known (see for instance
\cite{argyros2005ramsey}).

\begin{teo}
Let $\uf$ be an $\alpha$-uniform family over a set $M$. Then $\uf$
is lexicographically well ordered and its order type is
$\omega^\alpha$.
\end{teo}

In what follows, we consider an uniform family $\mathcal{F}$ on
$\mathbb N$ (or a final segment of $\mathbb N$) as topological space
by giving $\mathcal{F}$ the order topology respect to the
lexicographic order $<_{\lex}$.

Now we recall some known facts about the Cantor-Bendixson
derivative (CB derivative in short). Given a topological space $X$
and $A\subseteq X$, we let $A'$ be the set of all limit points
$x\in A$. Recursively, $A^{(0)}=A$,  $A^{(\alpha+1)}$ is
$(A^{(\alpha)})'$ and for $\alpha$ a limit ordinal, $A^{(\alpha)}$
is $\bigcap_{\beta<\alpha}A^{(\beta)}$. The least $\alpha$ such
that $A^{(\alpha)}=A^{(\alpha+1)}$ is called the CB index of $A$.
It is well know that $\omega^\alpha$ with the order topology has
CB index equal to $\alpha$.

An ordinal is said to be {\em indecomposable} if there are not
$\beta,\gamma <\alpha$ such that $\alpha=\beta+\gamma$. It is
known that $\alpha$ is indecomposable iff $\alpha=\omega^\beta$
for some $\beta$.

To get copies of uniform families we will use the following
theorem which follows from the results in \cite{Baumgart86}.

\begin{teo}
\label{copias_de_w_alfa} Let $\alpha<\omega_1$ be an
indecomposable ordinal and  $X\subseteq \omega^\alpha$. If
$X^{(\gamma)}\neq\emptyset$ for all $\gamma<\alpha$, then $X$ has
a subspace homeomorphic to $\omega^\alpha$.
\end{teo}

\section{CB derivatives of uniform families}

In this section we study the behavior of the CB derivative on
$\uf\upharpoonright M$, for $M\in\mathbb{N}^{[\infty]}$,  as a
subspace of $\mathcal{F}$. In particular, we will characterize the
limit points in $\mathcal{F}\upharpoonright M$.

\begin{lema}
\label{longitud} Let $\mathcal{F}$ be an $\alpha$-uniform family
on a final segment of $\mathbb{N}$ with $\alpha\geq\omega$ and $t\in\uf$.

\begin{itemize}
\item[$(i)$] If $\alpha=\omega$, then $|t|\geq \min(t)+1$.

\item[$(ii)$] If $\alpha>\omega$, then $|t|> \min(t)+2$.
\end{itemize}

In particular, $|t|\geq 2$ for all $t$ in an $\alpha$-uniform family
with $\alpha>1$ and $\min(t)\geq 1$.

\end{lema}

\begin{proof} Let $\mathcal F$ be an $\omega$-uniform family and $t\in
\mathcal F$. Let $n=\min(t)$, then $t/n\in\mathcal{F}_{\{n\}}$ and
$\mathcal{F}_{\{n\}}$ is $k$-uniform with $k\geq n$, therefore the
size of $t$ is at least $n+1$. The rest of the claim follows by
induction on $\alpha$.

\end{proof}

\begin{lema}
\label{convergencia_abajo} \label{converegencia_tn} Let
$\mathcal{F}$ be an uniform family on a final segment of
$\mathbb{N}$.

\begin{enumerate}

\item[$(i)$] Suppose $(s_i)_i$ is a sequence in $\mathcal F$  such
that $s_i\rightarrow s$ with $s\in \mathcal F$, then there exists
$k\in \mathbb{N}$ such that $\min(s)-1\leq \min(s_i)\leq \min(s)$
for all $i\geq k$. In particular, $s_i\leq_{\lex} s$ for all $i\geq
k$.

\item[$(ii)$] Suppose $(s_i)_i$ is a sequence in $\uf$ of the form
$s_i=u\cup\{p-1\}\cup v_i$ where $u\in\finite$,  $p\geq 1$,
$\max(u)<p-1<\min(v_i)$ and $\min(v_i)\rightarrow \infty$. Then
there is $m\in \mathbb N$ such that
\[
s_i\rightarrow u\cup\{p, p+1,\dots,p+m\}=u\cup t_p^{{\mathcal
F}_{u}}.
\]

\item[$(iii)$] Suppose $(s_i)_i$ is a sequence in $\mathcal F$
such that $s_i\rightarrow s\in\mathcal{F}$ and
$\min(s_i)=\min(s)-1=p-1$ for all $i$. Then $s=t_p^{\mathcal F}$ and
$s_i=\{p-1\}\cup v_i$ for some $v_i$ such that $p-1<\min(v_i)$ and
$\min(v_i)\rightarrow \infty$. Conversely,  if $s_i\rightarrow_i
t_p^{\mathcal F}$ and $s_i\neq t_p^{\mathcal F}$ for all $i$, then
eventually $\min(s_i)= p-1$.

\item[$(iv)$] Suppose $s_i\rightarrow s$ and
$\min(s_i)=\min(s)=n$ for all $i$. Then $s_i/n\rightarrow s/n$.

\item[$(v)$] Suppose $s, s_i\in \uf$ with $s\neq s_i$ for all $i$
and $s_i\rightarrow s$. Then there are $u, v_i\in \finite$ and
$p\in \mathbb N$ such that
\[
s=u\cup t_p^{\uf_u}
\]
and eventually
\[
s_i=u\cup \{p-1\}\cup v_i
\]
where  $\max(u)<p-1<\min (v_i)$ and $\min(v_i)\rightarrow \infty$.

\end{enumerate}
\end{lema}

\begin{proof} $(i)$ follows from the fact that $\mathcal{F}$ is a front and
the topology of $\mathcal F$ is the order topology given by
$<_{\lex}$ which is a well-order on $\mathcal{F}$. In particular,
convergence in $\mathcal F$ is from below.

To see $(ii)$, let $s=u\cup t_p^{{\mathcal F}_{u}}$  and $w\in
\mathcal F$ such that $w<_{\lex} s$.  It is clear that $s_i
<_{\lex}s$ for all $i$. We will show that eventually $w<_{\lex}
s_i$. The only interesting case is when $w=u\cup v$ with
$\max(u)<\min(v)$. If $\min(v)<p-1$, then clearly $w<_{\lex} s_i$
for all $i$. Suppose then that $\min(v)=p-1$.  As $s_i\in \uf$ and
$\mathcal F$ is a $\sqsubseteq$-antichain, then $u\cup\{p-1\}\not\in
\mathcal F$ and thus $|v|\geq 2$. Therefore, $w<_{\lex} s_i$ for all
large enough $i$.

For $(iii)$, notice that $s_i <_{\lex} t_p^{\mathcal F}\leq_{\lex}
s$ for all $i$. Thus $s=t_p^{\uf}$. Suppose $m$ is such that
$\min(v_i)< m$. Since $\uf$ is a front, pick $w_m\in \uf$ such that
$\{p-1, m\}\sqsubseteq w_m$. Then  $s_i<_{\lex} w_m <_{\lex} s$.
Hence there are only finitely many such $v_i$ and thus
$\min(v_i)\rightarrow \infty$.

To see $(v)$. By $(i)$ we assume that $s_i\leq_{\lex} s$ for all
$i$. If $\min(s_i)=\min(s)-1$ eventually, then apply $(iii)$ to get
the conclusion with $u=\emptyset$. If $\min(s_i)=\min(s)=n$, then by
$(iv)$, $s_i/n\rightarrow s/n$; by repeating this finitely many
times we get that $s=u\cup w$, $s_i=u\cup w_i$ with $\max(u)< w$,
$\max(u)<\min(w_i)$, $\min(w_i)=\min(w)-1$ and $w_i\rightarrow w$.
Since $w,w_i\in \uf_u$ and $\uf_u$ is uniform on ${\mathbb N}/u$,
then we apply $(iii)$ to finish the proof.

\end{proof}

\begin{rem}
{\em Let $\uf$ be an uniform family on $\mathbb N$. If ${\mathcal
F}\upharpoonright M$ is a closed subset of $\mathcal F$, then $M$
is a final segment of $\mathbb N$. In fact, let $n\in M$, we show
that $n+1\in M$. Since ${\mathcal F} \upharpoonright M$  is a
front on $M$, let $v_i\in{\mathcal F}\upharpoonright M$ such that
$\{n, i\}\sqsubseteq v_i$ for $i\in M/n$. Then $v_i\rightarrow
t_{n+1}^{\mathcal F}$, in particular $n+1\in M$.}
\end{rem}

Using the previous results, we are ready to characterize limit
points in uniform families.

\begin{prop}
\label{limites_familias_uniformes_M} Let $\mathcal{F}$ be an
$\alpha$-uniform family on a final segment $S$ of $\mathbb{N}$ with
$1<\alpha<\omega_1$, $M\in S^{^{[\infty]}}$ and
$t\in\mathcal{F}\upharpoonright M$ with $\min(t)>1$. Then,
$t\in(\mathcal{F}\upharpoonright M)'$ if, and only if, there is
$u\in  S^{^{[<\infty]}}$ and $p\in \mathbb N$ such that
\[
t=u\cup\{p,p+1,\cdots,p+m\}
\]
where $\max(u)<p-1$,  $p-1\in M$ and $m\geq 1$. Notice that $t=u\cup
t_p^{{\mathcal F}_u}$.
\end{prop}

\begin{proof} $(\Rightarrow)$ Let $t\in(\mathcal{F}\upharpoonright M)'$, by
Lemma \ref{converegencia_tn} we know that there is $u\in \finite$
and $p,m\in \mathbb N$ such that
\[
t=u\cup\{p,p+1,\cdots,p+m\} = u\cup t_p^{\uf_u}
\]
and $\max(u)<p-1$. Moreover, any sequence in $\uf\upharpoonright M$
converging to $t$ is eventually of the form $s_i=u\cup \{p-1\}\cup
v_i$ where $\max(u)<p-1<\min (v_i)$ and $\min(v_i)\rightarrow
\infty$. In particular, $p-1\in M$. It remains only  to show that
$m\geq 1$. Since $\{p-1\}\cup v_i\in \uf_u$, then $\uf_u$ is not
$1$-uniform, thus by Lemma \ref{longitud}, $t_p^{\uf_u}$ has size at
least 2, hence $m\geq 1$.

\medskip

$(\Leftarrow)$ Reciprocally, suppose $t=u\cup
t_p^{\mathcal{F}_u}\subseteq M$ for some $p\in M$ with
$\max(u)<p-1\in M$. Notice that ${\mathcal F_u}\upharpoonright M$ is
a $\beta$-uniform family on $M/u$ for some $\beta<\alpha$. Since
$t_p^{\uf_u}$ has size at least 2, then $\beta\geq 2$. As
$\uf_u\upharpoonright M$ is a front on $M/u$, there is $w_i\in
{\mathcal F_u}\upharpoonright M$ such that $\{p-1,i\}\sqsubseteq
w_i$ for each $i\in M/(p-1)$. Then by Lemma \ref{converegencia_tn}
we know that $u\cup w_i\rightarrow u\cup t_p^{\uf_u}$.

\end{proof}

\medskip

Proposition \ref{limites_familias_uniformes_M} gives a tool to
determine the topological type of a subspace
$\mathcal{F}\upharpoonright M$. Also, it allows to construct
subspaces $\mathcal{F}\upharpoonright M$ without copies of
$\mathcal{F}$. The following example shows that
$\mathcal{F}\upharpoonright M$ can be a discrete subspace of
$\mathcal{F}$.

\begin{ejems}
For the following examples we shall consider the Schreier barrier
$\mathcal{S}$ (defined in \S 2).
\begin{enumerate}
\item[$(i)$] Let $M\in\mathbb{N}^{^{[\infty]}}$ be the collection
of even numbers. Since in $M$ there are not  consecutive numbers,
then $\mathcal{S}\upharpoonright M$ is a discrete subspace of
$\mathcal{S}$.

\item[$(ii)$] Let $M=\{3k:k\in \mathbb{N}\}$ and
$N=\mathbb{N}\backslash M$. In this case, $N$ has consecutive
numbers but $\mathcal{S}\upharpoonright N$ is also discrete,
because $3q\notin N$ for all $q$.
\end{enumerate}
\end{ejems}

As we can see, given an uniform family $\mathcal{F}$ on $\mathbb
N$, its restrictions $\mathcal{F}\upharpoonright M$ can change
considerably its topological type. Nevertheless, for some sets $M$
the restriction conserves the topological type of $\mathcal{F}$.
The simplest example is when $M$ is a final segment of $\mathbb
N$, then $\mathcal{F}\upharpoonright M$ corresponds also to final
segment of $\mathcal{F}$, therefore $\mathcal{F}\upharpoonright M$
is closed in $\mathcal{F}$ and the subspace topology of
$\mathcal{F}\upharpoonright M$ is homeomorphic $\mathcal F$. But,
as we shall show in following sections, there are also non trivial
sets $M$ such that $\mathcal{F}\upharpoonright M$ contains a
topological copy of $\mathcal{F}$. To do this, we need to analyze
the CB derivatives of an uniform family.

Using the definition of $\mathcal{F}_{\{n\}}$, $\varoplus$, and
$<_{\lex}$, it is easy to verify the following result which we shall
use continuously to make proofs by induction.

\begin{lema}\label{F_n}
Let $\mathcal{F}\subseteq \finite$ and
$M\in\mathbb{N}^{^{[\infty]}}$. The following hold:
\begin{enumerate}
\item[$(i)$] $\mathcal{F}_{\{n\}}\upharpoonright
M=(\mathcal{F}\upharpoonright M)_{\{n\}}$, for $n\in M$,

\item[$(ii)$]
$\mathcal{F}_{\{n\}}=\bigcup_{m>n}(\mathcal{F}_{\{n\}})_{\{m\}}\varoplus
\{\{m\}\}$, for $n \in \mathbb{N}$,

\item[$(iii)$] $\mathcal{F}\upharpoonright M=\bigcup_{n\in
M}(\mathcal{F}\upharpoonright M)_{\{n\}}\varoplus \{\{n\}\}$.

\end{enumerate}
\end{lema}

\begin{lema}\label{derivada_entra}
Let $\mathcal{F}$ be an $\alpha$-uniform family on a final segment
$S$ of $\mathbb N$, $M$ an infinite subset of $S$, $u$ a finite set and
$0<\beta<\alpha$, then
$$
\big[(\mathcal{F}\upharpoonright M)_{u}\varoplus
\{u\}\big]^{^{(\beta)}}=\big[(\mathcal{F}\upharpoonright
M)_{u}\big]^{^{(\beta)}}\varoplus \{u\}.
$$
In particular, for $n\in \mathbb{N}$ we have
$$
\big[(\mathcal{F}\upharpoonright M)_{\{n\}}\varoplus
\{\{n\}\}\big]^{^{(\beta)}}=\big[(\mathcal{F}\upharpoonright
M)_{\{n\}}\big]^{^{(\beta)}}\varoplus \{\{n\}\}.
$$

\end{lema}

\begin{proof} By induction on $\beta$. The result its true for $\beta=1$ by
Lemma \ref{convergencia_abajo}.  Let us consider $\beta<\alpha$ and
let us suppose that the lemma is true for all $\gamma<\beta$.
\begin{enumerate}
\item[$(i)$] Suppose $\beta=\gamma+1$ and let
$t\in\big[(\mathcal{F}\upharpoonright M)_{u}\varoplus
\{u\}\big]^{^{(\gamma+1)}}$. Then there exists  $(t_i)_i$ in
\linebreak $ \big[(\mathcal{F}\upharpoonright M)_{u}\varoplus
\{u\}\big]^{^{(\gamma)}}$ such that $t_i\rightarrow t$. By the
inductive hypothesis,  \linebreak $(t_i)_i\in
\big[(\mathcal{F}\upharpoonright M)_{u}\big]^{^{(\gamma)}}\varoplus
\{u\}$. Thus, by Lemma \ref{convergencia_abajo}, we get that
\mbox{$t/u \in\big[(\mathcal{F}\upharpoonright
M)_{u}\big]^{^{(\gamma+1)}}$}. Hence $t=u\cup t/u \in
\big[(\mathcal{F}\upharpoonright
M)_{u}\big]^{^{(\gamma+1)}}\varoplus\{u\}$.

Reciprocally, let $t\in\big[(\mathcal{F}\upharpoonright
M)_{u}\big]^{^{(\gamma+1)}}\varoplus \{u\}$. Then $t/u \in
\big[(\mathcal{F}\upharpoonright M)_{u}\big]^{^{(\gamma+1)}}$. Thus there is $(t_i)_i \in
\big[(\mathcal{F}\upharpoonright M)_{u}\big]^{^{(\gamma)}}$
such that $t_i\rightarrow t/u$. Hence $t\in
\big[(\mathcal{F}\upharpoonright
M)_{\{n\}}\varoplus\{u\}\big]^{^{(\gamma+1)}}$, because
\[
u\cup t_i\in \big[(\mathcal{F}\upharpoonright
M)_{u}\big]^{^{(\gamma)}}\varoplus\{u\}=\big[(\mathcal{F}\upharpoonright
M)_{u}\varoplus\{u\}\big]^{^{(\gamma)}}.
\]

\item[(ii)] If $\beta$ is an ordinal limit, then

\notag\begin{align}
\hspace{-1cm}\big[(\mathcal{F}\upharpoonright M)_{u}\varoplus \{u\}\big]^{^{(\beta)}}&=\bigcap_{\lambda<\beta} \big[(\mathcal{F}\upharpoonright M)_{u}\varoplus \{u\}\big]^{^{(\lambda)}}\\
&=\bigcap_{\lambda<\beta} \Big(\big[(\mathcal{F}\upharpoonright M)_{u}\big]^{^{(\lambda)}}\varoplus \{u\}\Big)\\
&=\Big(\bigcap_{\lambda<\beta} \big[(\mathcal{F}\upharpoonright M)_{u}\big]^{^{(\lambda)}}\Big)\varoplus \{u\}\\
&=\big[(\mathcal{F}\upharpoonright M)_{u}\big]^{^{(\beta)}}\varoplus \{u\}.
\end{align}\qedhere
\end{enumerate}
\end{proof}

\begin{prop}\label{F_F_n}
Let $\mathcal{F}$ be an $\alpha$-uniform family on a final segment
of $\mathbb N$ with $2<\alpha<\omega_1$,
$M\in\mathbb{N}^{^{[\infty]}}$ and $0<\beta<\alpha$. If $t\in
(\mathcal{F}\upharpoonright M)^{^{(\beta)}}$ then one of the
following holds:
\begin{enumerate}
\item[$(i)$] $t/k\in ((\mathcal{F}\upharpoonright
M)_{\{k\}})^{^{(\beta)}}$, where $k=\min(t)$, or

\item[$(ii)$] $t=t_p^\mathcal{F}$, for some $p\in\mathbb{N}$ with
$p-1\in M$.
\end{enumerate}

Therefore

$$
(\mathcal{F}\upharpoonright M)^{^{(\beta)}}\subseteq\bigcup_{k\in
M}[(\mathcal{F}\upharpoonright
M)_{\{k\}}\varoplus\{\{k\}\}]^{^{(\beta)}} \cup
\{t_p^\mathcal{F}:\; t_p^\mathcal{F}\subseteq M\; \mbox{and}
\;p-1\in M\}.
$$
\end{prop}

\begin{proof} Note that the last equation is consequence of $(i)$, $(ii)$
and Lemma \ref{derivada_entra}. On the other hand, let $t\in
(\mathcal{F}\upharpoonright M)^{^{(\beta)}}$ and $k=\min(t)$. Then
$t_k^{\uf}\leq_{\lex} t<_{\lex} t_{k+1}^{\uf}$. There are two cases
to consider: (a) Suppose $t=t_k^\mathcal{F}$. Since $t$ is a limit
point, then by Lemma \ref{limites_familias_uniformes_M}, $k-1\in M$
and (ii) holds.

(b) Suppose $t_k^{\uf}<_{\lex} t$. Let
\[
U_k=\{s\in \mathcal{F}|M:\; t_k^{\mathcal F}<_{\lex} s<_{\lex}
t_{k+1}^{\mathcal F}\}.
\]
Then  $t\in U_k$ and $U_k$ is an open subset of
\makebox{$\mathcal{F}\upharpoonright M$. Thus $t\in
({U}_k)^{^{(\beta)}}\subseteq ((\mathcal{F}\upharpoonright
M)_{\{k\}}\varoplus \{\{k\}\})^{^{(\beta)}}$}
$=((\mathcal{F}\upharpoonright M)_{\{k\}})^{^{(\beta)}}\varoplus
\{\{k\}\}$. Thus (i) holds.

\end{proof}

\subsection{Finite CB derivative}

In what follows we present some results about the finite derivatives
$(\uf\upharpoonright M)^{^{(l)}}$, with $l<\omega$.

\begin{lema} \label{extension} Let $\mathcal{F}$ be an
$\alpha$-uniform family on $M$ with $\alpha\geq \omega$. There is
a sequence $(w_j)_j$  of finite sets with  $(min(w_j))_j$
increasing and an increasing sequence of integers $(k_j)_j$ such
that $\uf_{w_j}$ is $k_j$-uniform on $M/w_j$.
\end{lema}

\begin{proof} By induction on $\alpha$. For $\alpha=\omega$ the result
follows from the definition of a $\omega$-uniform family. If
$\alpha>\omega$, then  $\uf_{\{j\}}$ is $\beta_j$-uniform on $M/j$
with $\omega\leq \beta_j<\alpha$  for (eventually) all $j\in M$.
Using the inductive hypothesis, define recursively $k_j$ and $v_j$
for $j\in M$ such that $\uf_{\{j\}\cup v_j}$ is $k_j$-uniform on
$M/{v_j}$, $j<\min (v_j)$ and $(k_j)_j$ increasing. Take $w_j=
\{j\}\cup v_j$ with $j\in M$.
\end{proof}

\begin{prop}
\label{C-B_k-uniformes_coro2} Let $\uf$ be an $\alpha$-uniform
family on a final segment $S$  of $\mathbb N$ with \linebreak
$\alpha\geq 3$ and   $M \in S^{^{[\infty]}}$.  Suppose there is
$l\in\mathbb{N}$ with $1\leq l$ and $N \in \mathbb{N}^{^{[\infty]}}$
such that \linebreak $\{i,i+1,i+2,\dots,i+l\}\subseteq M$ for all
$i\in N$. Let $u\in \finite$ and $p>\max(u)+1$ be such that
$\uf_{u\cup\{p-1\}}$ is $\beta$-uniform with $l\leq\beta$. If $t\in
\uf$ is of the form
\[
t=u\cup\{p,p+1,\dots,p+m\}
\]
with $l\leq m$, then $t\in (\uf\upharpoonright M)^{^{(l)}}$.
\end{prop}

\begin{proof}  When $l=1$, the result follows from Proposition \ref{limites_familias_uniformes_M}, thus we assume $l\geq 2$. Let
$t$, $M$ and $N$ as in the hypothesis. We will define a sequence
$(s_i)_i$ in $(\uf\upharpoonright M)^{^{(l-1)}}$ converging to $t$.

We treat first the case $\beta<\omega$.  When $l=\beta$, take
$s_i= u\cup \{p-1\}\cup\{i+1, \cdots,i+l\}$ for $i\in N/p$. If
$l<\beta$, then for infinite many $i\in N$ there is a nonempty
finite set $w_i$ such that
\[
s_i=u\cup \{p-1\}\cup w_i\cup\{i+1, \cdots,i+l\}\in
\uf\upharpoonright M,
\]
$p-1<\min(w_i)$,  $\max(w_i)<i$ and $\min(w_i)\rightarrow \infty$.
This finishes the definition of the sequence $(s_i)_i$.  By a
straightforward inductive argument, we conclude that $s_i\in
(\uf\upharpoonright M)^{^{(l-1)}}$. By Lemma \ref{converegencia_tn},
$s_i\rightarrow t$  and thus $t\in(\uf\upharpoonright M)^{^{(l)}}$.

Now suppose $\beta\geq \omega$. By Lemma \ref{extension}, there are
sequences $(w_i)_i$ and $(k_i)_i$ such that
\mbox{$p<\min(w_i)\rightarrow \infty$}, $k_i>m$ and
$\uf_{u\cup\{p-1\}\cup w_i}$ is $k_i$-uniform. Then we construct the
sequence $(s_i)_i$ as before.
\end{proof}

For $k$-uniform families with $k\in\omega$ we have the following
proposition.

\begin{prop}
\label{C-B_k-uniformes_coro} Let $\uf$ be a $k$-uniform family on
a final segment of $\mathbb N$ with $3\leq k$.  Let
$l\in\mathbb{N}$ with $2\leq l<k$,  $M \in
\mathbb{N}^{^{[\infty]}}$ and \;$t \subseteq M$. If $t\in
(\uf\upharpoonright M)^{^{(l)}}$, then  there exist $N \in
\mathbb{N}^{^{[\infty]}}$ such that
$\{i,i+1,i+2,\dots,i+l\}\subseteq M$ for all $i\in N$ and
\[
t=u\cup\{p,p+1,\dots,p+m\}
\]
for some $u\in \finite$ with $\max(u)<p-1\in M$ and $l\leq m \leq
k-1$.
\end{prop}

\begin{proof}  Let $t\in (\uf\upharpoonright M)^{^{(l)}}$, then by
Proposition \ref{limites_familias_uniformes_M}
\[
t=u\cup\{p,p+1,\dots,p+m\}
\]
for some $u\in \finite$ with $\max(u)<p-1\in M$. Let $(s_i)_i$ in
$(\uf\upharpoonright M)^{^{(l-1)}}$ converging to $t$. By Lemma
\ref{converegencia_tn} we assume that each  $s_i$ is of the form
\[
s_i=u\cup\{p-1\}\cup v_i
\]
with $p-1<\min(v_i)$.

The proof is by induction on $l$. By the inductive hypothesis when
$l\geq 3$ and by Proposition \ref{limites_familias_uniformes_M} when
$l=2$, we conclude that there is an increasing sequence $(p_i)_i$
such that $p_i-1\in M$, $\{p_i, p_i+1,\cdots, p_i+m_i\}\subseteq
v_i$ and $l-1\leq m_i$. In particular, this says that
\mbox{$\{p_i-1, p_i, p_i+1,\cdots, p_i+l-1\}\subset M$} for all $i$.

Now we show that $l\leq m <|t|-1$.  In fact, $m=|t|-|u|-1=|v_i|\geq
m_i+1\geq l$.

\end{proof}

From the previous results we immediately get the following:

\begin{teo}\label{M_k}
Let $M \in \mathbb{N}^{^{[\infty]}}$ and $k> 2$. Then $M^{^{[k]}}$,
as a subspace of $\mathbb{N}^{^{[k]}}$, has \linebreak CB index $k$
if, and only if, there exists $p\in\mathbb{N}$ and $N \in
\mathbb{N}^{^{[\infty]}}$ such that\linebreak
$\{p-1,p,p+1,p+2,\dots,p+k-1\}\subseteq M$ and
$\{i,i+1,i+2,\dots,i+k-1\}\subseteq M$ for all $i\in N$.
\end{teo}\hfill{$\Box$}

The previous theorem  gives a characterization of those $M\in
\mathbb{N}^{^{[\infty]}}$ such that the CB index of
$\mathcal{F}=\mathbb{N}^{^{[k]}}$ and  $\mathcal{F}\upharpoonright
M$ are the same. However, this does not guarantee that
$\mathcal{F}\upharpoonright M$ contains a topological copy of
$\mathcal F$. To get this, we need that
$\{p-1,p,p+1,p+2,\dots,p+k-1\}\subseteq M$ for infinite many $p$.

The following example shows what we have said in the introduction
about Nash-Williams theorem.

\begin{ejem}
\label{nash-williams} Let $\uf$ be a $\alpha$-uniform family on
$\mathbb N$ with $\alpha\geq 2$. Let ${\mathcal B}=\uf^{(1)}$ and
$M$ be an infinite set. We will show that
$(\uf\upharpoonright M) \setminus {\mathcal B}\neq\emptyset$. In
particular, this says that if $M$ is homogeneous for the partition
given by $\mathcal B$, then  $(\uf\upharpoonright M)\cap \uf^{(1)}=\emptyset$ and thus $\uf\upharpoonright M$ is a
discrete subset of $\uf$.

Suppose first that $\alpha\geq \omega$. By Lemma \ref{extension},
applied to $\uf\upharpoonright M$, there is $u\subset M$ finite such
that $\uf_u\upharpoonright M$ is $k$-uniform for some $2\leq
k<\omega$. Let $w\subset M$ and $p,q\in M$ such that $\max(w)<p<q-1$
and $|w\cup\{p,q\}|=k$. Then $t=u\cup w\cup\{p,q\}\in
\uf\upharpoonright M$ and $t\not\in \mathcal B$ (by Proposition
\ref{limites_familias_uniformes_M}). If $\alpha<\omega$, we can
argue analogously to find $t$.

\end{ejem}

\section{$\uf$-adequate sets}

Let $\uf$ be an $\alpha$-uniform family on a final segment $S$ of
$\mathbb N$ with $\alpha\geq 2$. In this section we introduce the
notion of a $\uf$-adequate set $M$ and later  we will show that for those sets
$\uf\upharpoonright M$ has the same CB index as $\uf$.

Let $M \in S^{^{[\infty]}}$, we define  by recursion a subset
$M(\uf)$ of $M$ and the notion of a $\uf$-adequate set.

\begin{itemize}
\item[$(i)$] If $\alpha=2$, then $M(\uf)$ is the set of all $n\in
M$ such that $t_{n+1}^{\uf}\subset M$. And $M$ is said to be
$\uf$-adequate, if $M(\uf)$ is not empty.

\item[$(ii)$] If $\alpha=\beta+1$, then
\[
M(\uf)=\{n\in M:\; t_{n+1}^{\uf}\subset M,  \mbox{$M/n$ is
$\uf_{\{n\}}$-adequate and $(M/n)(\uf_{\{n\}})$ is infinite}\}.
\]
And $M$ is said to be $\uf$-adequate, if $M(\uf)$ is not empty.

\item[$(iii)$] If $\alpha$ is limit, then $M(\uf)=M$. Let
$(\alpha_n)_n$ be the increasing sequence of ordinals as in the
definition of a $\alpha$-uniform family. We say that $M$ is
$\uf$-adequate, if for all $n$ there is a non empty finite set
$v\subset M$ such that $\uf_v$ is $\gamma$-uniform for some
$\gamma\geq \alpha_n$ and $M/v$ is $\uf_v$-adequate.
\end{itemize}

\begin{ejem}
\label{ejemplo} If $\uf=\mathbb{N}^{[2]}$, then an infinite set is
$\uf$-adequate when it contains three consecutive integers. In
general, for $\uf=\mathbb{N}^{[k]}$, a set is $\uf$-adequate if it
contains $\{n,n+1, \cdots, n+k\}$ for some $n$ and infinite many
intervals of length $k$.

Let us say that an infinite set $M$ is $\omega$-adequate, if it
contains arbitrarily long intervals of consecutive integers.
Suppose $\uf$ is $\omega$-uniform on $\mathbb N$. Then $M$ is
$\uf$-adequate iff $M$ is $\omega$-adequate.

Now suppose that $\uf$ is $(\omega+1)$-uniform on $\mathbb N$. Let
$P$ be a $\omega$-adequate set. For a fixed $k\in \mathbb N$, let
$M=P\cup\{k\}\cup t_{k+1}^{\uf}$. Then $M$ is $\uf$-adequate. In
fact, notice that $k\in M(\uf)$ because $M/k$ is $\omega$-adequate
and $\uf_{\{k\}}$ is $\omega$-uniform.
\end{ejem}

The next lemma says that, in the definition of a $\uf$-adequate
set for $\alpha$ limit, we could have required that the ordinals
$\gamma$ are successor.

\begin{lema}\label{adecuado-limite}
Let $\uf$ be an $\alpha$-uniform family on a final segment of
$\mathbb N$ with $\alpha$ a limit ordinal. If $M$ is an
$\uf$-adequate set, then there is  a sequence of ordinals
$\beta_n<\alpha$ and  finite sets $u_n\subset M$ such that $M/u_n$
is $\uf_{u_n}$-adequate, $\uf_{u_n}$ is $(\beta_n+1)$-uniform on
$M/u_n$, $\alpha=sup\{\beta_n:\; n\in \mathbb{N}\}$.
\end{lema}

\begin{proof} By induction. The result holds for $\alpha=\omega$ by the
definition of an $\omega$-uniform family. Let $\alpha>\omega$ be a
limit ordinal. Let $(\alpha_n)_n$ converging to $\alpha$ as in the
definition of an $\alpha$-uniform family.  Fix sequences
$(\gamma_n)_n$ and $(v_n)_n$ as in the definition of $\uf$-adequate
set.  Since $(\alpha_n)_n$ is increasing, we assume that
$\gamma_n>\alpha_n$. If there are infinitely many $n$ such that
$\gamma_n$ is a successor ordinal, then we are done. Otherwise,
assume that $\gamma_n$ is a limit ordinal for all $n$. Apply the
inductive hypothesis to $\uf_{v_n}$ and $M/v_n$ to  get sequences of
ordinals $\beta^n_k$ converging to $\gamma_n$ and finite sets
$v^n_k\subset M$ such that $v_n\sqsubset v^n_k$, $M/v^n_k$ is
$\uf_{v^n_k}$-adequate and $\uf_{v^n_k}$ is $(\beta^n_k+1)$-uniform.
Now pick for each $n$ an integer $k_n$ such that
$\beta^n_{k_n}>\alpha_n$. Take $u_n= v^n_{k_n}$ and
$\beta_n=\beta^n_{k_n}$.
\end{proof}

We are going to present a method to construct $\uf$-adequate sets.

We need to introduce a notation. If $\mathcal B$ is a collection
of finite sets, then $\overline{\mathcal{B}}^{\; \sqsubseteq}$
denotes the collection of all finite sets $t$ such that
$t\sqsubseteq s$ for some $s\in \mathcal B$.

It is easy to show by induction on $\alpha$ that if $\mathcal{F}$ is
$\alpha$-uniform with $\alpha\geq\omega$, then there exist $s\in
\finite$ such that $\mathcal{F}_s$ is $\omega$-uniform on
$\mathbb{N}/s$. Thus the following definition is non trivial.

\begin{defi}
Let $\mathcal{F}$ be an $\alpha$-uniform family with
$\alpha\geq\omega$, we define the set $\mathcal{A}_\mathcal{F}$ as

$$
\mathcal{A}_\mathcal{F}=\{s\in\overline{\mathcal{F}}^{\;
\sqsubseteq}: \mathcal{F}_s \text{ is $\omega$-uniform on
}\mathbb{N}/s\}.
$$
\end{defi}

The set $\mathcal{A}_\mathcal{F}$ has the following properties:

\begin{enumerate}
\item $\mathcal{A}_\mathcal{F}$ is infinite, if
$\alpha\neq\omega$,

\item $\mathcal{A}_\mathcal{F}$ is a front on $M$ (If
$\mathcal{F}$ is uniform on $M\in\mathbb{N}^{^{[\infty]}})$,

\item $\overline{\mathcal{A}_\mathcal{F}}^{\;\sqsubseteq}$ is a
well founded tree.
\end{enumerate}

From $\mathcal{A}_\mathcal{F}$ we define a $\mathcal{F}$-adequate
tree and then a $\mathcal{F}$-adequate set of natural numbers.

\begin{defi}
Let $\mathcal{F}$ be an $\alpha$-uniform family with
$\alpha\geq\omega$. We will say that a non  empty subset $T$ of
$\overline{\mathcal{A}_\mathcal{F}}^{\; \sqsubseteq}$ is a
\textbf{$\mathcal{F}$-tree}, if the following conditions hold

\begin{enumerate}
\item[$(i)$] If $t\in T$ and $s\sqsubseteq t$, then $s\in T$,

\item[$(ii)$] $Ter(T)\subseteq\mathcal{A}_\mathcal{F}$,

\item[$(iii)$] $\{n\in\mathbb{N}: n>t \text{ and }t\cup\{n\} \in
T\}$ is infinite, for all $t\in T\backslash Ter(T)$,
\end{enumerate}
where $Ter(T)$ denotes the set of terminal nodes of $T$.
\end{defi}

We remark that for an $\alpha$-uniform family $\uf$ on a set $M$ with $\alpha>\omega$,
$\mathcal{A}_\mathcal{F}$ is a front on
$M$, and thus $\overline{\mathcal{A}_\mathcal{F}}^{\;\sqsubseteq}$ is well founded \cite{argyros2005ramsey}. Thus, each
$\mathcal{F}$-tree is also well founded.

\begin{defi}
Given $\mathcal{F}$ an $\alpha$-uniform family with
$\alpha>\omega$ and $T$ a $\mathcal{F}$-tree, we define
$E(T)\in\mathbb{N}^{^{[\infty]}}$ as

$$
E(T)=\bigcup_{\begin{matrix}
                  s\cup\{n\}\in T \\
                  s<n
                \end{matrix}}\{n\}\cup t_{n+1}^{\mathcal{F}_s}.
$$

In other words,
$$
\emptyset\neq\{x_0,x_1,x_2,\dots x_{k-1},x_k\}\in T
\Leftrightarrow\{x_k\}\cup
t_{x_k+1}^{\mathcal{F}_{\{x_0,x_1,x_2,\dots x_{k-1}\}}}\subseteq
E(T).
$$
\end{defi}

The following result is easy to verify. 

\begin{lema}\label{A_f}
Let $\mathcal{F}$ be an $\alpha$-uniform family over a final segment of $\mathbb N$ with
$\alpha>\omega$ and $n\in\mathbb{N}$. Then,

\begin{enumerate}
\item $(\mathcal{A}_\mathcal{F})_{\{n\}}=\mathcal{A}_{\mathcal{F}_{\{n\}}}$,

\item If $T$ is a $\mathcal{F}$-tree, then $T_{\{n\}}$ is a
$\mathcal{F}_{\{n\}}$-tree for all $n$ such that $\{n\}\in T$,

\item $E(T_{\{n\}})\subseteq E(T)$  for all $n$ such that
$\{n\}\in T$.
\end{enumerate}
\end{lema}

\begin{prop}
\label{F-tree-adecuado} Let $\uf$ be an $\alpha$-uniform family
over a final segment of $\mathbb N$ with $\alpha>\omega$. If $T$
is a $\uf$-tree, then $E(T)$ is $\uf$-adequate.
\end{prop}

\begin{proof} By induction on $\alpha$. Let us fix a $\uf$-tree $T$ and let
$M=E(T)$.   We will show that  $M$ is $\uf$-adequate and moreover
that it is infinite.

\begin{itemize}
\item[(i)] Suppose $\alpha=\omega+1$. It is easy to verify that
$n\in M$ for all $n$ such that $\{n\}\in T$. Recall that by Lemma
\ref{longitud}, the size of $t_{n+1}^{\uf}$ is increasing with $n$.
Thus, $M$ contains arbitrarily long intervals of consecutive
integers and by Example \ref{ejemplo}, $M$ is $\uf_{\{n\}}$ adequate
for all $n$.

\item[(ii)] If $\alpha =\beta+1$, we will show that $M(\uf)$
contains all $n$ such that $\{n\}\in T$. Fix such an $n$.  Then
$t_{n+1}^{\uf}\subset M$. Let $M_n$ be $E(T_{\{n\}})$. Since
$T_{\{n\}}$ is a $\uf_{\{n\}}$-tree, by the inductive hypothesis,
$M_n$ is $\uf_{\{n\}}$-adequate and $M_n(\uf_{\{n\}})$ is infinite.
As \mbox{$M_n(\uf_{\{n\}})\subset M_n \subset M/n$}, then $M/n$ is
$\uf_{\{n\}}$-adequate. Thus $n\in M(\uf)$.

\item[(iii)] Finally, suppose $\alpha$ is a limit ordinal. Then
$T_{\{n\}}$ is a $\uf_{\{n\}}$-tree for each $n$ such that $\{n\}\in
T$. Since $\uf_{\{n\}}$ is $\alpha_n$-uniform, then $E(T_{\{n\}})$
is  $\uf_{\{n\}}$-adequate. Since $E(T_{\{n\}})\subseteq E(T)$, then
$E(T)$ is also $\uf_{\{n\}}$-adequate. As this holds for infinite
many $n$'s, then $E(T)$ is $\uf$-adequate.\end{itemize}\end{proof}

\begin{ejem}
Let $\uf$ be  a $(\omega+1)$-uniform family on $\mathbb N$. It is
easy to construct an infinite set $P$ containing arbitrarily long
intervals of consecutive natural numbers and such that
$t_n^{\uf}\not\subset P$ for all $n$. As in Example \ref{ejemplo},
fix $k\in \mathbb N$ and  let $M=P\cup\{k\}\cup t_{k+1}^{\uf}$. Then
$M$ is $\uf$-adequate and it is not of the form $E(T)$ for any
$\uf$-tree $T$.
\end{ejem}

\section{Topological copies of $\uf$ inside $\mathcal{F}\upharpoonright M$}

The following theorem is one of the main results of this paper. It
justifies the introduction of $\uf$-adequate sets.

\begin{teo}
\label{rango de adecuados} Let $\uf$ be an $\alpha$-uniform family
on a final segment $S$ of $\mathbb N$ with $\alpha\geq 2$ and  $M$
a $\uf$-adequate set. Then the CB index of $\uf\upharpoonright M$
is $\alpha$.

\end{teo}

\begin{proof}  Since $\uf$ is homeomorphic to $\omega^\alpha$, then the CB
index of $\uf\upharpoonright M$ is at most $\alpha$.

We first show by induction on $\beta\geq 1$ that if $\uf$ is
$(\beta+1)$-uniform, $M$ is $\uf$-adequate and $n\in M(\uf)$, then
\[
t_{n+1}^{\mathcal{F}}\in (\mathcal{F}\upharpoonright
M)^{^{(\beta)}}.
\]

\begin{enumerate}

\item[$(i)$] If $\beta=1$, then  $t_{n+1}^{\uf}=\{n+1,n+2\}\subset
M$. From Proposition \ref{limites_familias_uniformes_M},
$t_{n+1}^{\uf} \in(\mathcal{F}\upharpoonright M)^{(1)}$.

\item[$(ii)$] Suppose $\beta=\gamma+1$. Since $M/n$ is
$\uf_{\{n\}}$-adequate and $(M/n)(\uf_{\{n\}})$ is infinite, there
is an increasing sequence $k_i\in (M/n)(\uf_{\{n\}})$. Then by the
inductive hypothesis, $t_{k_i+1}^{\uf_{\{n\}}}\in
(\mathcal{F}_{\{n\}}\upharpoonright M)^{^{(\gamma)}}$. By Lemma
\ref{derivada_entra} we have
\[
s_i=\{n\}\cup t_{k_i+1}^{\uf_{\{n\}}}\in
(\mathcal{F}\upharpoonright M)^{^{(\gamma)}}.
\]
By Lemma \ref{converegencia_tn}, $s_i\rightarrow t_{n+1}^{\uf}$.
Thus $t_{n+1}^{\uf}\in (\mathcal{F}\upharpoonright
M)^{^{(\gamma+1)}} $ and we are done.

\item[$(iii)$] Suppose $\beta$ is a limit ordinal. Let
$\beta_m\uparrow\beta$ as in the definition of a $\beta$-uniform
family. Since $M/n$ is $\uf_{\{n\}}$-adequate, then there is a
sequence of finite sets $u_m\subset M/n$ and ordinals
$\gamma_m\geq\beta_m$ such that $\mathcal{G}_m=\uf_{\{n\}\cup u_m}$
is $\gamma_m$-uniform on $M/u_m$ and $M/u_m$ is
$\mathcal{G}_m$-adequate. By Lemma \ref{adecuado-limite}, we assume
that each $\gamma_n$ is a successor ordinal. Let $k_m\in
M(\mathcal{G}_m)$. Then by the inductive hypothesis
$t_{k_m+1}^{\mathcal{G}_m}\in (\mathcal{G}_m\upharpoonright
M)^{^{(\beta_m)}}$. By Lemma \ref{derivada_entra} we have
\[
s_m=\{n\}\cup u_m\cup t_{k_m+1}^{\mathcal{G}_m}\in
(\mathcal{F}\upharpoonright M)^{^{(\beta_m)}}.
\]
By Lemma \ref{converegencia_tn}, $s_m\rightarrow t_{n+1}^{\uf}$.
Thus $t_{n+1}^{\uf}\in (\mathcal{F}\upharpoonright M)^{^{(\beta)}} $
and we are done.

\end{enumerate}

The proof of the theorem is by induction on $\alpha$. It remains
only to consider the case when $\alpha$ is a limit ordinal. Let
$(\alpha_k)_k$ be an increasing sequence of ordinals converging to
$\alpha$ as in the definition of a $\alpha$-uniform family. Since
$M$ is $\uf$-adequate, then for all $k$ there is a finite set
$v_k\subset M$ such that $M/v_k$ is $\uf_{v_k}$-adequate and
$\uf_{v_k}$ is $\gamma_k$-uniform with $\gamma_k\geq \alpha_k$. By
the inductive hypothesis, the CB index of $\uf_{v_k}\upharpoonright
M/v_k$ is $\gamma_k$ and therefore (by Lemma \ref{derivada_entra})
the CB index of $\uf\upharpoonright M$ is larger than $\gamma_k$ for
all $k$. Thus this last index is $\alpha$.
\end{proof}

For $\alpha=\omega$ we have a more complete result than Theorem
\ref{rango de adecuados} as follows.

\begin{teo}\label{S_M_adecuado_teo}
Let $\mathcal{F}$ be a $\omega$-uniform family on a final segment
of $\mathbb N$ and $M \in \mathbb{N}^{^{[\infty]}}$. Then,
$\mathcal{F}\upharpoonright M$ has CB index $\omega$ if, and only
if, $M$ is $\uf$-adequate.
\end{teo}

\begin{proof} The {\em if} part follows from Theorem \ref{rango de adecuados}. For
the other direction we will use the characterization of
$\uf$-adequate sets given in Example \ref{ejemplo}.

Let $\mathcal{F}$ be a $\omega$-uniform family on $S$ and $(m_k)_k$
be an strictly increasing sequence in $\mathbb{N}$ such that
$\mathcal{F}_{\{k\}}$ is $m_k$-uniform on $S/k$ for all
$k\in\mathbb{N}$. Suppose  $\mathcal{F}\upharpoonright M$ has CB
index $\omega$. Then, given $n\in\mathbb{N}$ there exists $t\in
\big(\mathcal{F}\upharpoonright M\big)^{^{(n)}}$ and a sequence
$(t_i)_i$ in $\big(\mathcal{F}\upharpoonright~M\big)^{^{(n-1)}}$
such that $t_i\uparrow t$. Let $k_i=\min(t_i)$, by Proposition
\ref{F_F_n}, for all $i\in\mathbb{N}$,  $t_i/k_i\in
\big((M/k_i)^{^{[m_{k_i}]}}\big)^{^{(n-1)}}$ or\;
$t_i=t_{k_i}^\mathcal{F}$ with \makebox{$k_i-1\in M$}. Since
$(t_i)_i$ is convergent, then eventually $t_i\neq
t_{k_i}^\mathcal{F}$. Therefore, by Proposition
\ref{C-B_k-uniformes_coro}, we can suppose that each $t_i/k_i$ has
the form $t_i/k_i= u_i\cup \{p_i,p_i+1,\dots,p_i+n-1\}$ with
$p_i-1\in M$ for each $i\in \mathbb{N}$. Hence,
$\{p_i-1,p_i,p_i+1,\dots,p_i+n-1\}\subseteq M$ for all
$i\in\mathbb{N}$, which implies $M$ is $\uf$-adequate.
\end{proof}

\begin{coro} \label{S_M_adecuado_coro}

Let $\mathcal{F}$ be a $\omega$-uniform family and $M \in
\mathbb{N}^{^{[\infty]}}$. Then, $\mathcal{F}\upharpoonright M$
has a topological copy of $\mathcal{F}$ if, and only if, $M$ is
$\uf$-adequate.\end{coro}

\begin{proof} Let $\mathcal{F}$ be a $\omega$-uniform family and $M \in
\mathbb{N}^{^{[\infty]}}$. If  $\mathcal{F}\upharpoonright M$
contains a topological copy of $\mathcal{F}$, then
$\mathcal{F}\upharpoonright M$ has CB index $\omega$ and therefore
by Theorem \ref{S_M_adecuado_teo} $M$ is $\uf$-adequate.
Reciprocally, if $M$ is an $\uf$-adequate set, then by Theorem
\ref{S_M_adecuado_teo} $\mathcal{F}\upharpoonright M$ has CB index
$\omega$, and by Theorem \ref{copias_de_w_alfa}
$\mathcal{F}\upharpoonright M$ has a topological copy of
$\mathcal{F}$.
\end{proof}

Finally, we present a result about the restriction to a set of the
form $E(T)$ for $T$ a $\uf$-tree.

\begin{teo}
Let $\mathcal{F}$ be an $\alpha$-uniform family with
$\alpha>\omega$ indecomposable. If $T$ is a $\mathcal{F}$-tree,
then $\mathcal{F}\upharpoonright E(T)$ contains a topological copy
of $\mathcal{F}$.
\end{teo}

\begin{proof} Let $\mathcal{F}$, $\alpha$ and $T$ be as in the hypothesis.
Then by Proposition \ref{F-tree-adecuado}, we know that $E(T)$ is
$\uf$-adequate. Hence by Theorem \ref{rango de adecuados},
$\mathcal{F}\upharpoonright E(T)$ has CB index $\alpha$, and by
Theorem \ref{copias_de_w_alfa}, $\mathcal{F}\upharpoonright E(T)$
has a topological copy of $\mathcal{F}$.

\end{proof}

\bibliographystyle{plain}

\end{document}